\documentclass{amsart}
\usepackage{amsmath}
\usepackage{amssymb}
\usepackage{fancybox}
\usepackage{eucal}
\usepackage{amsthm}
\usepackage{amscd}
\usepackage{hyperref}
\usepackage{wasysym}
\usepackage{url}

\usepackage{amssymb,}
\usepackage[]{amsmath, amsthm, amsfonts,graphicx, amscd,}
\usepackage[all,cmtip]{xy}
\input amssym.def \input amssym

\newtheorem {thm}{Theorem}[section]

\newtheorem{prop}[thm]{Proposition}

\newtheorem{ques}[thm]{Question}

\theoremstyle{remark}
\newtheorem{rem}[thm]{Remark}
\newtheorem{np*}{Non-Proof}

\theoremstyle{definition}

\newcommand{\pd}[2]{\frac{\partial #1}{\partial #2}}

\def\Ind{\setbox0=\hbox{$x$}\kern\wd0\hbox to 0pt{\hss$\mid$\hss} \lower.9\ht0\hbox to 0pt{\hss$\smile$\hss}\kern\wd0}
\def\Notind{\setbox0=\hbox{$x$}\kern\wd0\hbox to 0pt{\mathchardef \nn=12854\hss$\nn$\kern1.4\wd0\hss}\hbox to 0pt{\hss$\mid$\hss}\lower.9\ht0 \hbox to 0pt{\hss$\smile$\hss}\kern\wd0}

\newcommand{\m}{\mathbb }

\begin{document}

\title{Painlev\'e equations, vector fields, and ranks in differential fields}
\author{James Freitag} 
\address{James Freitag\\ Department of Mathematics, Statistics, and Computer Science\\
University of Illinois at Chicago\\
322 Science and Engineering Offices (M/C 249)\\
851 S. Morgan Street\\
Chicago, IL 60607-7045}
\email{freitagj@gmail.com}

\begin{abstract}
  Model theoretic ranks of solutions to Painlev\'e equations are calculated, and the type of the generic solution of the second Painlev\'e equation is shown to be disintegrated, strengthening a theorem of Nagloo. A question of Hrushovski and Scanlon regarding Lascar rank and Morley rank in differential fields is answered using planar vector fields. 
\end{abstract}

\thanks{Thanks to Ronnie Nagloo and Anand Pillay for useful discussion around the Painlev\'e equations, as well as many enlightening lectures on the subject over the past several years. Thanks to Dave Marker and Matthias Aschenbrenner for useful conversations on ranks in differential fields and the Painlev\'e equations.}

\maketitle

\section{Introduction} 
The purpose of this short note is to answer a question of Hrushovski and Scanlon, while at the same time explaining why that question was still open before the appearance of this note. In \cite{hrushovski1999lascar}, Hrushovski and Scanlon gave a series of examples of definable sets in differentially closed fields with Lascar rank not equal to Morley rank, answering a question of Poizat. Around the same time, as noted in \cite{hrushovski1999lascar}, Marker and Pillay proved that for differential varieties of algebraic dimension two over the constants, Lascar rank equals Morley rank (by \emph{algebraic dimension} we mean the Kolchin polynomial of a generic point on the variety - in this context, the algebraic dimension is given by the order of the equation being considered). The minimal algebraic dimension of the examples of \cite{hrushovski1999lascar} is five. It is easy to see that differential varieties of algebraic dimension one have Lascar rank and Morley rank one. In light of Marker and Pillay's argument, it seemed plausible that all differential varieties of algebraic dimension two might have Lascar rank equal to Morley rank. So, naturally Hrushovski and Scanlon \cite[Question 2.9]{hrushovski1999lascar} ask if there is a theorem explaining the gap between two and five. 

Nagloo and Pillay \cite{nagloo2011algebraic} claimed that the total space of the second Painlev\'e equation is an algebraic dimension three set witnessing that Lacar rank is not equal to Morley rank; their paper \cite{nagloo2011algebraic} used deep model theoretic tools to strengthen results of the Umemura and Watanabe (cited in detail below), proving results about orthogonality of various fibers of the Painlev\'e families. However, in the case of the calculation of the ranks, only a slight reinterpretation of the results of Umemura and Watanabe was required. Following this, the present author and Moosa \cite{JRJou} proved that for arbitrary algebraic dimension two differential varieties, Lascar rank equals Morley rank. These developments would have entirely closed the question posed by Hrushovski and Scanlon; the gap would have been explained by \cite[Theorem 6.1]{JRJou} and counterexamples coming from the second Painlev\'e equation. 

However, we will show that the claims about rank of \cite{nagloo2011algebraic} were not correct; for any of the families of Painlev\'e equations, Lascar rank equals Morley rank. The problem seems to have been largely caused by differences in language. Our rank calculations follow from the work of Umemura and Watanabe, and we will give careful and specific citations for each of the Painlev\'e families. So, the following question gives the remaining open content of the question of Hrushovski and Scanlon:

\begin{ques} \label{qq}
	Is there a differential variety with algebraic dimension three (or four) with Lascar rank not equal to Morley rank?
\end{ques} 

We answer Question \ref{qq} with a new example of an order three differential variety with Lascar rank not equal to Morley rank. This automatically also gives an order four example (by taking the product with any order one differential variety). Our example is based certain rational vector fields on the affine plane, and is related to recent work generalizing results of Rosenlicht to the case of equations with nonconstant coefficients and the Poincar\'e problem \cite{JRosen}. 

Our rank calculations for the Painlev\'e equations also allow one to give a series of examples which answer the following natural question; such examples seem to not have appeared in the literature previously:

\begin{ques} \label{degcont} Is there a family of differential varieties $X \rightarrow Y$ in which for some $n \in \m N $, $\{ y \in Y \, | \, d_M (X_y) =n \} $ is not definable ($d_M(-)$ denotes the Morley degree)?  
\end{ques}

The rank calculations also allow for an easy generalization of the arguments of \cite{nagloo2015geometric} proving the geometric triviality of the fibers of the second Painlev\'e family.

Morley degree one is a natural notion of irreducibility for definable sets, but the Kolchin topology already comes equipped with its own notion of irreducibility:  a Kolchin-closed set $X \subset \m A^n$ over $K$ is irreducible if the collection of differential polynomials with coefficients in $K$ which vanish on $X$ forms a prime differential ideal. 

\begin{ques} \label{ritt} Is there a family of differential varieties $X \rightarrow Y$ in which $\{ y \in Y \, | \, X_y \text{ is irreducible} \} $ is not definable? 
\end{ques} 
Question \ref{ritt} is equivalent to the \emph{Ritt problem}, an important open problem which has received considerable attention \cite{golubitsky2009generalized}. Questions \ref{ritt} and \ref{degcont} are of a similar flavor, but the interaction of model-theoretic ranks with the Kolchin topology is somewhat enigmatic; for instance, \cite{freitag2012model} gives an example of a definable set whose Kolchin closure has higher Lascar and Morley rank than the original set. 

\subsection{Organization} 
Section \ref{Correc} is devoted to explaining the details of the rank calculations for the Painlev\'e equations. Section \ref{examp} is devoted to answering Question 2.9 of \cite{hrushovski1999lascar}.

\section{The Painlev\'e equations} \label{Correc}
\subsection{Notation} 
Let $K$ be a differential field with derivation $\delta$ which contains an element $t$ such that $\delta t =1$. The following subsections will be devoted to rank calculations of each of the families of Painlev\'e equations. 

\subsection{The Painlev\'e two} 
The differential algebraic information in this section comes from \cite{umemura1997solutions}, whose notation we also follow. The second Painlev\'e family of differential equations is given by 
$$P_{II} (\alpha) : \, \, y'' = 2y^3 +ty + \alpha$$ where $\alpha$ ranges over the constants. For $\alpha=-\frac{1}{2}$, Umemura and Watanabe \cite[see 2.7-2.9 on pages 169--170]{umemura1997solutions} show that if $K_1$ is a differential field extension of $K$ and if $y_1$ is a solution to $P_ {II} ( -\frac{1}{2}) $ such that the transcendence degree of $K_1 \langle y_1 \rangle $ over $K_1$ is one, then $y_1$ satisfies the Ricatti equation: 
$$y_1 ' = y_1 ^2 +\frac{1}{2} t.$$ In model theoretic terms, this implies that $$\left\{y \, | \, y'' = 2y^3 +ty -\frac{1}{2} , \,  y' \neq y^2 + \frac{1}{2} t \right\}$$ is strongly minimal, while $P_{II} ( -\frac{1}{2})$ is of Morley rank one and Morley degree two. The differential varieties $P_{II} ( \alpha)$ for $\alpha \in \frac{1}{2} + \m Z$ are all isomorphic via Backlund transformations, so the same analysis applies to $P_{II} ( \alpha )$ for $\alpha \in \frac{1}{2} +\m Z $. Note however, that the degree of the exceptional subvariety changes with the application of the Backlund transformations (this can be seen by direct calculations, though the fact that the degree can not be bounded over all $\alpha \in \frac{1}{2} +\m Z $ also follows by a standard compactness argument). 

This contradicts the remarks in subsection 3.7 of \cite{nagloo2011algebraic}, where it is claimed that the Morley rank (and Lascar rank) of $P_{II}(\alpha)$ for $\alpha \in \frac{1}{2} + \m Z$ is two (in particular Fact 3.22 of \cite{nagloo2011algebraic} is incorrect). Parts of the subsequent discussion of the subsection depend on this fact, and it is easy to see that this leaves the order three case Question 2.9 of \cite{hrushovski1999lascar} still open. The Painlev\'e II family witnesses the non-definability of Morley degree, rather than Morley rank. In the coming subsections, we will show that this is the case for each of the Painlev\'e families. 

For Painlev\'e II, we can also expand a result of Nagloo \cite{nagloo2015geometric}: 

\begin{prop} The definable set 
$$X= \left\{ y \, \vert \, y'' = 2y^3 +ty -\frac{1}{2} , \,  y' \neq y^2 + \frac{1}{2} t \right\}= P_{II} \left( -\frac{1}{2} \right) \setminus \left\{ y \, \vert \, y' = y^2 +\frac{1}{2} t \right\}$$ 
is strongly minimal and geometrically disintegrated.
\end{prop} 
\begin{proof} As established in the paragraphs above, the strong minimality of this set is a reinterpretation of the results of Umemura and Watanabe \cite[169--170]{umemura1997solutions}. With this in place, we will establish the disintegratedness of this definable set via the argument of Nagloo \cite[Proposition 3.3]{nagloo2015geometric}. 
	
By the strong minimality of the above set $X$, the type of a generic solution to $P_{II} ( -\frac{1}{2})$ is of Morley rank one. The equivalence relation of nonorthogonality refines transcendence degree, so the type of a generic solution to $P_{II} ( -\frac{1}{2})$ is orthogonal to the constants. By a result of Hrushovski and Sokolovic \cite{HrSo}, any locally modular nondisintegrated strongly minimal set in differentially closed fields is nonorthogonal to the Manin kernel of a simple abelian variety. From this, it follows that a strongly minimal set $X$ is disintegrated if for any generic $x,y \in X$, if $y \in K \langle x \rangle ^{alg}$, then $y \in K \langle x \rangle $ (see \cite[2.7]{nagloo2015geometric} for a proof). The remaining portion of the proof follows \cite[Proposition 3.3]{nagloo2015geometric} almost verbatim; it can be verified that the strong minimality of $P_{II}(\alpha)$ is used in only one place in the proof of Proposition 3.3 of \cite{nagloo2015geometric}. Namely, in Claim 1 of the proof of Proposition 3.3, strong minimality is only used to show that the polynomial $F$ (defined in \cite[Claim 1 in the proof of Proposition 3.3]{nagloo2015geometric}) cannot divide its derivative. The same applies in our case by strong minimality of $X$. The rest of the argument proceeds identically to \cite[Proposition 3.3]{nagloo2015geometric}. 
	\end{proof} 
	
Because the Backlund transformations give definable bijections between the fibers of $P_{II} (\alpha)$ for $\alpha \in \frac{1}{2} + \m Z$, it is the case that for each such fiber, the generic type of each fiber of the the second Painlev\'e equation is geometrically trivial. 

\subsection{Painlev\'e three} 
The differential algebraic information in this section comes from \cite{umemura1998solutions}, whose notation we also follow. For the purposes of determining the Morley rank and degree of the fibers of the third Painlev\'e family, it is sufficient consider the following system of equations, which we denote by $S( \bar v)$, 
\begin{eqnarray*} 
t q ' & = &  2 q^2 p -q^2 -v_1 q +t  \\ 
t p' &= & -2 qp^2 + 2 qp -v_1 p + \frac{1}{2} (v_1 +v_2)
\end{eqnarray*}
which can be obtained from the third Painlev\'e family via the transformation given in the introduction of \cite{umemura1998solutions}, where $\bar v = (v_1, v_2) \in \m C^2$. 
Define $$W_1 = \{ \bar v \in \m C^2 \, | \, v_1+v_2 \in 2 \m Z \text{ or } v_1 -v_2 \in 2 \m Z \}$$ and $$D_1 = \{ \bar v \in \m Z^2 \, | \, v_1 +v_2 \in 2 \m Z\}.$$ 
Theorem 1.2 (iii) \cite{umemura1998solutions} implies that for $\bar v$ not in $W_1$ or $D_1$, $S(\bar v)$ is strongly minimal. For $\bar v \in W_1$, Lemma 3.1 of \cite{umemura1998solutions} implies that $S( \bar v)$ has Morley rank one and Morley degree two.  For $\bar v \in D_1$, Lemma 3.2 of \cite{umemura1998solutions} implies that $S( \bar v)$ has Morley rank one and Morley degree three. To see these latter two facts from the statements of Lemmas 3.1 and 3.2, note that the fibers of the family related by an affine transformation in the group generated by
$$s_1 (\bar v) = (v_2, v_1) , \, s_2 (\bar v) = (-v_2 , -v_1), \, s_3 (\bar v ) = (v_2+1, v_1-1), \, s_4 ( \bar v ) = (-v_2+1, -v_1+1)$$ are isomorphic.


\subsection{Painlev\'e four} 
The differential algebraic information in this section comes from \cite{umemura1997solutions}, whose notation we also follow. Again let $K$ be a differential field with derivation $\delta$ which contains an element $t$ such that $\delta t =1$. Assume that the field of constants of $K$ is the field of complex numbers. The fourth Painlev\'e family of equations is given by $$y ''  = \frac{1}{2y}(y')^2 + \frac{3}{2}y^3 + 4tq^2 + 2 (t^2 - \alpha) y + \frac{\beta }{y},$$ where $ \alpha , \beta $ range over the constants. We denote the solution set to the previous equation by $P_{IV} ( \alpha , \beta )$. 

Let $S_{IV} ( v_1, v_2, v_3)$ be the solution set to the following system of differential equations: 

\begin{eqnarray*}
q ' &=& 2pq -q^2 -2tq + 2 ( v_1 -v_2) \\
p' &=& 2pq-p^2 +2tp+ 2(v_1-v_3)
\end{eqnarray*}

where $\bar v = (v_1, v_2, v_3 )$ are constants such that $\bar v \in V := \{ \bar v = \, | \, v_1+v_2+v_3=0 \} $. Then when $ \alpha = 3v_3 +1 , \beta = -2 (v_2 -v_1 )^2$ the elements $q$ such that there is a $p$ so that $(q,p ) \in (S_{IV} ( v_1, v_2, v_3)$ is precisely $P_{IV} (\alpha , \beta )$. 

Define the following affine transformations: 
\begin{eqnarray*} 
s_1 (v_1, v_2, v_3) & =& (v_2, v_1, v_3)\\
s_2 (v_1, v_2, v_3) &=& (v_3,v_2,v_1)	\\
 t_{-}  (v_1, v_2, v_3) &=&  (v_1, v_2, v_3)  + \frac{1}{3} (-1,-1,2) \\ 
 s_0 &=& t_- ^{-1} s_1 s_2 s_1 t_-
  \end{eqnarray*} 
  Let $H$ be the subgroup generated by $s_0,s_1,s_2$. 
  
  Let $\Gamma $ be the subset of $\m C^3$ such that 
  \begin{eqnarray*} 
  Re( v_2 -v_1 ) & \geq & 0 \\
   Re( v_1 -v_3 ) & \geq & 0 \\
    Re( v_3 -v_2 +1 ) & \geq & 0 \\
    Im ( v_2 -v_1 ) & \geq &  0 \, \text{  if } Re (v_2-v_1)=0 \\
     Im ( v_1 -v_3 ) & \geq & 0 \,  \text{  if } Re( v_1 -v_3 ) =0 \\
     Im ( v_3 -v_2 ) & \geq &  0 \, \text{  if } Re( v_3 -v_2 +1 )=0
  \end{eqnarray*}
  
The set $\Gamma$ is a fundamental region of $V$ for the group $H$. For parameters $\bar v, \bar w$ which are in the same orbit under $H$, the sets $S_{IV} (\bar v) $ and $S_{IV} (\bar w) $ are isomorphic, so, to analyze the fourth Painlev\'e family, it will only be necessary to analyze those $\bar v \in \Gamma$. 

Define $$W = \{ \bar v \in V \, | \, v_1 -v_2 \in \m Z \text{ or } v_3 -v_2 \in \m Z \text{ or } v_1 -v_3  \in \m Z \}.$$ Corollaries 3.5 and 3.9 of \cite{umemura1997solutions} imply that $S_{IV} (\bar v)$ satisfies Condition J (has no differential subvarieties except for finite sets of points) when $\bar v \in \Gamma \setminus W$. It is easy to see that Condition $J$ is equivalent to strong minimality \cite[see the appendix for a discussion]{nagloo2011algebraic}. 

Define $$D = \{ \bar v \in V \, | \, v_1 -v_2 \in \m Z \text{ and } v_3 -v_2 \in \m Z \text{ and } v_1 -v_3  \in \m Z \}.$$ Then noting that $D$ is the orbit of the origin under $H$, Lemma 3.11 \cite{umemura1997solutions} implies that for $\bar v \in D$, $S_{IV}(\bar v)$ has two irreducible order one differential subvarieties over any differential field $K$ extending $\m C (t)$. So, $S_{IV}(\bar v)$ has Morley rank one and Morley degree 3. If $\bar v \in W \setminus D$, then Lemma 3.10 \cite{umemura1997solutions} implies that $S_{IV} (\bar v)$ has one irreducible order one differential subvariety, and so $S_{IV} (\bar v)$ has Morley rank one and Morley degree 2. 


\subsection{Painlev\'e five} 
The differential algebraic information in this section comes from \cite{watanabe1995solutions}, whose notation we also follow. The fifth Painlev\'e family is equivalent to the following system of equations \begin{eqnarray*} t Q' & = & 2 Q (Q-1)^2 P + (3v_1+ v_2 ) Q^2 -(t+4v_1)Q + v_1 -v_2 \\ tP' &=& (-3Q^2 +4Q -1) P^2 -2(3v_1 +v_2 ) Q P + (t+4v_1) P - (v_3-v_1)(v_4-v_1)
\end{eqnarray*} 
where $\bar v = (v_1, v_2, v_3, v_4) \in \m C^4$ lies in the hyperplane $V$ defined by $\sum v_i = 0$. 
Setting $q = \frac{Q}{Q-1 } $ and $p = -(Q-1)^2 P + (v_3-v_1 )(Q-1)$ then $p$ and $q$ satisfy 
\begin{eqnarray*} t q' & = & 2 q^2 p - 2 qp + t q^2 -tq + ( v_1-v_2-v_3+v_4) q + v_2 -v_1 \\ 
t p' &=& -2 qp^2 + p^2 -2tpq +tp -(v_1-v_2-v_3+v_4) p + (v_3-v_1)t
\end{eqnarray*} 

and solutions to this system are birational with the solutions to our earlier system. The properties we study are not sensitive up to birationality, so we will work with this second system, which we denote by $S ( \bar v)$.

Let $$W = \left\{ \bar v \in V \, | \, v_1 -v_2 \in \m Z \text{ or } v_1 -v_3 \in \m Z \text{ or } v_1 -v_4 \in \m Z \text{ or } v_2 -v_3 \in \m Z \text{ or } v_2 -v_4 \in \m Z \text{ or } v_3 -v_4 \in \m Z \right\}.$$

Corollary 2.6 of  \cite{watanabe1995solutions} implies that for $\bar v \notin W$, $S( \bar v )$ is strongly minimal. In particular, for generic parameters, $S( \bar v)$ is strongly minimal. Lemmas 3.1-3.4 of \cite{watanabe1995solutions} imply that for $\bar v \in W$ the Morley rank of $S (\bar v)$ has Morley rank one and Morley degree between two and four (the specific loci with a given degree can be easily deduced from the cited lemmas and noting that a group of affine transformations specified in \cite{watanabe1995solutions} acts on the family of equations).





\subsection{Painlev\'e six} 
Let $R$ be the collection of $24$ vectors of the following form: 
$$ (\pm 1 , \pm 1, 0 , 0) , (\pm 1 ,, 0  \pm 1,  0) , (\pm 1 ,0,0, \pm 1) , (0, \pm 1 , \pm 1,  0) , (0, \pm 1 ,0 , \pm 1) , (0,0, \pm 1 , \pm 1) .$$ Let $\langle \bar v , \bar w \rangle = v_1  \bar w_1 +  v_2  \bar w_2 +  v_3  \bar w_3 +  v_4  \bar w_4$ denote the usual inner product on $\m C^4$. For $\alpha \in R$ and $k \in \m Z$, define $$H_{\alpha, k } = \{ \bar v \in \m C^4 \, | \, \langle \bar v , \alpha \rangle = k \}.$$ 

Define $M$ to be the union of all $H_ { \alpha , k}$ for $ \alpha \in R$ and $k \in \m Z$. Let $P$ be the union of all intersections of the form $H_ { \alpha , k} \cap  H_ { \beta , l}$ such that $\alpha , \beta $ are linearly independent. Let $L$ be the union of all intersections of the form $H_ { \alpha , k} \cap  H_ { \beta , l} \cap H_{\gamma , m } $ for $ \alpha , \beta , \gamma \in R$ are linearly independent and $k,l,m \in \m Z$. Let $D$ be the union of all intersections of the form $H_ { \alpha , k} \cap  H_ { \beta , l} \cap H_{\gamma , m } \cap H_{\delta , n} $ for $ \alpha , \beta , \gamma, \delta  \in R$ are linearly independent and $k,l,m, n \in \m Z$. 

Theorem 2.1 (v) of \cite{watanabe1998birational} implies that for $\bar v \notin M$, $S (\bar v)$ is strongly minimal. Propositions 4.1, 4.4, and 4.9 imply that for $\bar v \in M \setminus P$, $S(\bar v)$ is Morley rank one and Morley degree two unless $v_1 -v_2 \in \frac{1}{2} + \m Z$ and $v_3-v_4 \in \m Z$, in which case $S( \bar v )$ has Morley degree four. Propositions 4.2 and 4.5 imply that for $\bar v \in P \setminus L$, $S(\bar v)$ is Morley rank one and Morley degree three. Propositions 4.3 and 4.6 imply that for $\bar v \in L \setminus D$, $S(\bar v)$ is Morley rank one and Morley degree three. Proposition 4.4 implies that for $\bar v \in L \setminus D$, $S(\bar v)$ is Morley rank one and Morley degree four. Proposition 4.7 implies that for $\bar v \in D$, $S( \bar v)$ is Morley rank one and Morley degree five. As in the previous sections, deriving the various above mentioned facts depends in each of the cases on the action of a group of affine transformations on the family of equations.

\section{Lascar rank and Morley rank at order three} \label{examp} 

Consider the differential variety $X_c$ given by: \begin{eqnarray} \label{dx} x' = cy+y-c \\
\label{dy} y' = \frac{y(y-1)}{x},
\end{eqnarray}
where $c$ is an arbitrary constant. A generic solution $(x,y)$ over $\m Q (c)$ generates a differential field extension of transcendence degree two. We aim to analyse the possible forking extensions of $(x,y)$. 

Let $F$ be an algebraically closed differential field extending $\m Q (c)$. If $x \in F$, then by equation \ref{dx}, as long as $c \neq -1$, we must have $y \in F$. Similarly, if $y \in F$, then $y' \in F$, and since $(x,y)$ is generic over $\m Q (c)$, we must have $y \neq 0, 1$ and $x \neq 0$, so by equation \ref{dy}, we must have $x \in F$. Thus, whenever $c \neq -1$, the forking extensions of $(x,y)$ over $F$ such that $F \langle x,y \rangle $ has transcedence degree one over $F$ must have the property that $(x,y)$ are interalgebraic over $F$. We will analyze the possible interalgebraicities of $x$ and $y$ over $F$ next. 

By Seidenberg's embedding theorem \cite{seidenberg1958abstract}, we may embed $\m Q(c) \langle x,y \rangle $ into the field of meromorphic functions on some domain $U \subset \m C$, thus regarding $x,y$ as functions of a complex variable $t$. Then 
\begin{eqnarray*}  \frac{dx}{dt} = cy+y-c \\
	 \frac{dy}{dt} = \frac{y(y-1)}{x}.
\end{eqnarray*}

Thus 
 \begin{eqnarray} \label{imp} \frac{dy}{dx} = \frac{y(y-1) }{x(cy+y-c)}.
\end{eqnarray}

The reader may verify that solutions to \ref{imp} satisfy the relation 
 \begin{eqnarray} \label{imp1} c_1 + \log (x) = c \log (y) + \log ( 1-y).
 \end{eqnarray}
 where $c_1$ is an arbitrary constant. Thus, when $0 \neq c \in \m Q$ and $e^{c_1}$ is transcendental and so $$e^{c_1} x = y^c (1-y)$$ witnesses that $(x,y)$ forks with $e^{c_1}$ over $\m Q$. Thus, when $c \in \m Q$, the Lascar rank $X_c$ is two. 
 
 When $c \not \in \m Q$, and $x,y$ are meromorphic functions satisfying equation \ref{imp1}, $x$ and $y$ are not both interalgebraic and transcendental over any field extension $F$ of $\m Q (c)$. So, by our above analysis, any forking extension of the type of $(x,y)$ over $\m Q$ is algebraic, so $X_c$ has Lascar rank one. By \cite[Theorem 6.1]{JRJou}, $X_c$ also has Morley rank one. Thus our family of varieties $\{X_c \, | \, c \in \m C\}$ witnesses the non-definability of Morley rank. Let $X \rightarrow \m C$ be the total family of differential varieties in which the fiber above $c \in \m C$ is $X_c$. When $c \in \m Q$, $X_c$ is an irreducible Kolchin closed set of Morley rank two. For any $c_1, c_2 \in \m Q$, $RM(X_{c_1} \cap X_{c_2}) <2$ since $X_{c_1} \cap X_{c_2}$ is a proper Kolchin closed subset of $X_{c_1}$. Thus, $RM(X) =3$. On the other hand, the Lascar rank of the generic type of the family is two, since the fiber of $X$ above any generic constant is strongly minimal. Any specialization of the generic type has transcendence degree at most two, and thus Lascar rank at most two. Thus, the Lascar rank of $X$ is two. By combining \cite[Theorem 6.1]{JRJou} and our example, \cite[Question 2.9]{hrushovski1999lascar} is now settled entirely. 
 
\begin{rem} Assume for simplicity in the following remark that $c \in \m N$. The general case in which $c \in \m Q$ is similar, but requires some slightly messier expressions attained after suitable exponentiation. 
	
 In the example above, for $c \in \m N$, $(x,y) \mapsto \frac{y^c(y-1)}{x}$ gives a definable map from the solution set of equation \ref{imp} to the constants. In the language of \cite[particularly, see the appendix]{JRosen}, $\frac{y^c(y-1)}{x}$ is a rational first integral. The reader can verify that $$\frac{y(y-1) }{x(cy+y-c)}= \frac{-\pd{ \left(\frac{y^c(y-1)}{x}\right)}{x}}{\pd{\left(\frac{y^c(y-1)}{x}\right)}{y}   },$$ and so by Theorem A.6 of \cite{JRosen}, the solution set is not weakly orthogonal to the constants. In the case that $c $ is irrational, no expression of $\frac{y(y-1) }{x(cy+y-c)}$ as a quotient of partial derivatives is possible, and one can deduce the fact that the solution set to \ref{imp} is orthogonal to the constants in this case from Theorem 1.2 of \cite{JRosen}. 
	\end{rem}

\bibliography{/Users/freitagj/Dropbox/Research}{}
\bibliographystyle{plain}

\end{document}